\documentclass[12pt]{article}
\begin{document}
\title{Modelling the Zero Point Field}
\author{B.G. Sidharth\\
Centre for Applicable Mathematics \& Computer Sciences\\
B.M. Birla Science Centre, Adarsh Nagar, Hyderabad - 500 063 (India)}
\date{}
\maketitle
\begin{abstract}
We construct a model of the Zero Point Field in terms of an infinite collection
of oscillators. This has relevance because of the recent identification of
Dark energy with such a Zero Point Field.
\end{abstract}
\section{Introduction}
The existence of the Zero Point Field (ZPF) was realised even by Max Planck. 
In its later version it has come to be known as Quantum Vaccuum, and
lately it has been possible to identify it with the mysterious dark
energy. The Zero Point Field arises because of the fact that though
classically, an oscillator in the ground state has zero energy, Quantum
Mechanically, owing to the uncertainity principle, there is an energy
fluctuation about the zero energy level.\\
Interestingly there have been two schools of thought. Zero Point Field, 
according to Quantum theoriests is a secondary effect arising from the
already present oscilaltors. However according to what has come to be
known as Stochastic Electrodynamics, the ZPF is primary, and infact
Quantum Mechanics can be deduced therefrom.  This is a chicken and egg
situation.\\
Recently the ZPF has come into focus once again because of the 
observed ever expanding, accelerating feature of the universe. In other
words there is a large scale repulsion represented by a cosmological constant 
$\Lambda$ once invented and then
rejected by Einstein, and this cosmic repulsion could be attributed to a
mysterious all pervading dark energy, which can be identified with the
ZPF: this could be the mechanism which drives the cosmic
expansion and acceleration (Cf.ref.\cite{r1}) for a 
review).
\section{Realising the ZPF}
Two of the earliest realisations of the ZPF were in the form of the Lamb shift
and the Casimir effect.\\
In the case of the Lamb shift, as is well known, the motion of an orbiting
electron gets affected by the background ZPF. Effectively there is an
additioinal field, over and above that of the nucleus. This additional
potential, as is well known is given by
 \cite{r2}
$$\Delta V (\vec r) = \frac{1}{2} \langle (\Delta r)^2 \rangle \nabla^2 V
(\vec r)$$
The additional energy
$$\Delta E = \langle \Delta V (\vec r) \rangle$$
contributes to the observed Lamb shift which is $\sim 1000 mc/sec$.\\
The essential idea of the Casimir effect is that the interaction between
the ZPF and matter leads to macroscopic consequences. For example if we
consider two parallel metallic plates in a conducting box, then we should
have a Casimir force given by
 \cite{r3}
$$F = \frac{-\pi^2}{240} \frac{\hbar cA}{l^4}$$
where $A$ is the area of the plates and $l$ is the distance between them.
More generally, the Casimir force is a result of the
boundedness or deviation from a Euclidean topology of or in the Quantum
Vaccuum. These Casimir forces have been experimentally demonstrated
\cite{r4,r5,r6,r7}.\\
Returning to the ZPF
as the ubiquitous dark energy, we observe that \cite{r8}, a fluctuating electromagnetic
field can be modelled as an infinite collection of independent harmonic oscillators.
Quantum Mechanically, the ground state of the harmonic oscillators is described by
$$\psi (x) = \left(\frac{m\omega}{\pi \hbar}\right)^{1/4} 
e^{-(m\omega/2\hbar)x^2}$$
which exhibits the probability for the oscillator to fluctuate, mostly in the region
given by
$$\Delta x \sim (\hbar /m\omega )^{1/2}$$
An infinite collection of such oscillators can be modelled by
$$\psi (\xi_1,\xi_2,\cdots ) = const. \exp [-(\xi^2_1 + \xi^2_2 + \cdots )],$$
which gives the probability amplitude for an electromagnetic field configuration
$B(x,y,z), \xi_1$, etc. being the Fourier coefficients. Finally, as a consequence there
is a fluctuating magnetic field given by
\begin{equation}
B = \frac{\sqrt{\hbar c}}{l^2}\label{eD}
\end{equation}
where $l$ is the extent over which the fluctuation is measured. Further
these fluctuations typically take place within the time $\tau$, a typical
elementary particle Compton time (Cf.ref.[1]). This begs the question whether such
ubiquotous fields could be tapped for terrestrial applications or otherwise.\\
We now invoke the well known result from macroscopic physics that the
current in a coil is given by
\begin{equation}
\imath = \frac{NBA}{R\Delta t}\label{eE}
\end{equation}
where $N$ is the numer of turns of the coil, $A$ is its area and $R$ the
resistance.\\
Introducing (\ref{eD}) into (\ref{eE}) we deduce that a coil in the ZPF
would have a fluctuating electric current given by
\begin{equation}
\imath \approx \frac{NA}{R} \cdot \frac{e}{l^2\tau}\label{eF}
\end{equation}
In principle it should be possible to harness the current (\ref{eF}).

\end{document}